\documentclass[12pt,reqno]{amsart}
\usepackage{amsmath, amssymb, amsthm}
\usepackage{euscript}
\usepackage{verbatim}
\usepackage{graphicx}

\newcommand{\Mbar}{M^{*}}

\DeclareMathOperator{\psibar}{\overline{\psi}}
\DeclareMathOperator{\grad}{grad}
\DeclareMathOperator{\inj}{inj}
\DeclareMathOperator{\diam}{diam}
\newcommand{\Ebar}{\overline{E}}
\newcommand{\Ehat}{\hat{E}}
\newcommand{\E}{E}
\newcommand{\bR}{\mathbb{R}}
\newcommand{\bZ}{\mathbb{Z}}

\newcommand{\bS}{\mathbb{S}}

\newcommand{\cC}{ {\mathcal C} }
\newcommand{\cB}{ {\mathcal B} }
\newcommand{\cTC}{ {\mathcal TC} }

\newcommand{\Minf}{M(\infty)}

\newcommand{\MlK}{M \backslash K}
\newcommand{\gbar}{\overline{g}}

\newcommand{\grd}{\mathring{g}}
\newcommand{\grdi}{ (\mathring{g}_i) }
\newcommand{\tilth}{\widetilde{\theta}}

\newcommand{\tilrp}{\widetilde{r}_p}
\newcommand{\tilrq}{\widetilde{r}_q}
\newcommand{\tilpp}{\widetilde{p}'}
\newcommand{\tilqp}{\widetilde{q}'}
\newcommand{\rarrow}{\rightarrow}
\newcommand{\inn}[2]{\left\langle  #1 \;, \; #2 \right\rangle}
\newtheorem{theorem}{Theorem}
\newtheorem{prop}[theorem]{Proposition}
\newtheorem{lemma}[theorem]{Lemma}

\begin{document}
\title[H\"older Compactification]{H\"older Compactification for some \\ manifolds with pinched negative \\ curvature near infinity}
\author[Bahuaud]{Eric Bahuaud}
\author[Marsh]{Tracey Marsh}
\address{University of Washington, Seattle, Washington}
\date{\today}

\begin{abstract} We consider a complete noncompact Riemannian manifold $M$ and give conditions on a compact submanifold $K \subset M$ so that the outward normal exponential map off of the boundary of $K$ is a diffeomorphism onto $\MlK$.  We use this to compactify $M$ and show that pinched negative sectional curvature outside $K$ implies $M$ has a compactification with a well defined H\"older structure independent of $K$.  The H\"older constant depends on the ratio of the curvature pinching.  This extends and generalizes a 1985 result of Anderson and Schoen.
\end{abstract}
\maketitle
The Poincar\'e model of hyperbolic space has a natural geometric compactification -- one can compactify by adding the sphere at infinity.  Taking this to be a model case for other simply connected manifolds of negative curvature leads to a classical construction made precise in \cite{Eberlein-ONeill}:  Let $M$ be a Cartan-Hadamard manifold, that is, a complete, simply connected Riemannian manifold with nonpositive sectional curvature.  Define an equivalence relation on the set of geodesic rays parametrized by arc length by saying that geodesic rays $\sigma$ and $\tau$ are \textit{asymptotic} if $d_M( \sigma(t), \tau(t))$ remains bounded as $t \rarrow +\infty$.  Here $d_M$ is the distance function on $M$ induced by the metric $g$.  We define $\Minf$ to be the set of all equivalence classes of this relation; this is the \textit{geometric boundary at infinity}. \label{page:asymp}

In 1985, Michael Anderson and Richard Schoen proved that given a Cartan-Hadamard manifold $M$ with pinched sectional curvatures like
\begin{eqnarray*}
-\infty < -b^2 \leq \sec( M ) \leq -a^2 < 0, 
\end{eqnarray*} 
where $a$ and $b$ are positive constants, then geometric boundary at infinity has a $C^{a/b}$ structure \cite{Anderson-Schoen}.  Motivated by this result we investigate to what extent the simply connected hypothesis may be relaxed when compactifying the manifold and what resulting regularity may be obtained for the compactified manifold with boundary.  In particular, let $M$ be a complete, noncompact Riemannian $(n+1)$-manifold.  Define an \textit{essential subset} \label{page:es} $K$ of $M$ to be a compact $(n+1)$-dimensional Riemannian submanifold with boundary, such that $Y := \partial K$ is a smooth hypersurface that is convex with respect to the outward pointing unit normal and such that $\exp\colon N^{+}Y \rarrow \overline{\MlK}$ is a diffeomorphism, where $N^{+}Y \approx Y \times [0, \infty)$ is the outward normal ray bundle of $Y$.  The main result of this paper is:
\begin{theorem} Let $(M,g)$ be a complete, noncompact Riemannian $(n+1)$-manifold.  Suppose that there exists  $K \subset M$, a compact $(n+1)$-dimensional Riemannian submanifold with boundary that satisfies:
\begin{enumerate}
\item $K$ is totally convex in $M$, i.e. if $p, q \in K$ and $\gamma\colon[0,1] \rarrow M$ is any geodesic with $\gamma(0) = p$ and $\gamma(1) = q$, then $\gamma( [0,1]) \subset K$.
\item $\MlK$ satisfies the following curvature assumption: 
\begin{eqnarray}
\label{curvasmp}
-\infty < -b^2 \leq \sec( \MlK ) \leq -a^2 < 0.
\end{eqnarray}
\end{enumerate}
Then $K$ is an essential subset of $M$ and $\Mbar := M \cup \Minf$ is a geometric compactification of $M$ as a topological manifold with boundary.  The boundary is homeomorphic to $\partial K$.  Further, $\Mbar$ is endowed with the structure of a $C^{a/b}$ manifold with boundary, independent of the choice of $K$.
\label{thm:main}
\end{theorem}

Since any point $x$ in a Cartan-Hadamard manifold $M$ may be regarded as a pole, any small closed ball about $x$ is easily seen to be an essential subset for $M$.  Therefore Theorem \ref{thm:main} generalizes and strengthens the Anderson-Schoen result, for it allows for much greater variety in the topology of $M$; essential subsets relax the stringent hypothesis of simple connectedness used in the Anderson-Schoen paper. In addition, the result here proves the regularity of the \textit{entire} compactification $\Mbar = M \cup \Minf$.  The Anderson-Schoen theorem only proves the regularity for the boundary at infinity $\Minf$.

The outline of this paper is as follows.  In Section \ref{sec:notation} we outline our notation and explain our comparison theorems.  In Section \ref{sec:tconvex} we provide a condition for an essential subset.  In Section \ref{sec:tracey} we describe the compactification of $M$ as a topological manifold, and then in Section \ref{sec:compreg} we set up the necessary estimates to show the compactified manifold has a well defined $C^{a/b}$ structure.

We would like to thank our advisor Jack Lee for suggesting this problem to us and for his constant guidance.  We would also like to thank Robin Graham for suggestions on an early draft of this paper.

\section{Notation and basic estimates}
\label{sec:notation}

In this section we outline our notation and provide the estimates that will be used in the subsequent comparison geometry.  Throughout this paper $M$ denotes a complete, connected, noncompact Riemannian $(n+1)$-manifold with metric $g$.  The letter $K$ will always denote an essential subset and $Y:=\partial K$ will denote the smooth hypersurface boundary.  Throughout this paper we assume the curvature assumption \eqref{curvasmp} and write $\alpha = a/b$.  There is no loss in generality in assuming that the pinching constants in \eqref{curvasmp} satisfy $a \leq 1 \leq b$, for the ratio of maximum to minimum sectional curvature, $a/b$, is invariant under a homothety of the metric.  Further, we follow the curvature sign conventions given in \cite{Lee} or \cite{Petersen}: for orthonormal vectors $X,Z$, the sectional curvature of the plane they span is $ \sec(X,Z) = R(X,Z,Z,X)$, where $R$ is the Riemannian curvature 4-tensor.

We trivialize the normal ray bundle with respect to the outward unit normal for $Y$ as $N^{+}Y \approx Y \times [0, \infty)$.  The exponential map restricted to $N^{+}Y$ is written $\E$.  For $p \in Y$, the notation $\gamma_p$ denotes the geodesic normal to $Y$ emanating outwards from $p$.  We call a geodesic ray \textit{untrapped} if it eventually escapes any compact set.

Following Petersen \cite{Petersen} we reserve the term \textit{geodesic segment} for a distance minimizing geodesic curve between two points.

It is easy to verify that $E\colon Y \times [0, \infty) \rarrow M$ is a local diffeomorphism at every point of $Y \times \{0\}$.  Therefore by compactness of $Y$ we may obtain an $\epsilon > 0$ and a one-sided collar neighbourhood $T$ of $Y$ so that $E\colon Y\times [0, \epsilon) \rarrow T$ is a diffeomorphism.  Let $r\colon T \rarrow \bR $ denote the distance to $Y$.  Adapting the proof of the classical Gauss lemma (see \cite{Lee} or \cite{Klingenberg} for a precise statement of the classical form) we obtain a decomposition of the metric as $g = dr^2 + g_Y(y,r)$.  Further, if we choose any coordinates $\{y^{\beta}\}$ on an open set $U \subset Y$ we may get coordinates on $T$ by extending $y^{\beta}$ to be constant along the integral curves of $\grad r$, and then $(y^{\beta}, r)$ form coordinates on $E(U \times [0, \epsilon)) \subset T$.  We will refer to such coordinates as \textit{Fermi coordinates for $Y$}, and in Section \ref{sec:tracey} we will see that total convexity implies Fermi coordinates for $Y$ exist on neighbourhoods of the form $U \times [0, \infty)$.

We use Latin indices to index directions in $M$ and consequently these indices range from $0$ to $n$.  We use Greek indices to index the directions along $Y$ which range from $1$ to $n$, and a zero or `$r$' to index the direction normal to $Y$.  

We now consider the second fundamental form of $Y$ in $M$ and we fix signs by taking our definition as $h(X,Z) = g( \nabla_X Z, -\partial_r)$, where $\nabla$ is the connection in $M$, and $X, Z$ are vector fields on $Y$ extended arbitrarily to vector fields on $M$.  Note that this definition uses the \textit{inward} pointing unit normal.  Given this convention, we say $Y$ is \textit{convex} (respectively \textit{strictly convex}) with respect to the outward unit normal if the scalar second fundamental form is positive semidefinite (respectively positive definite).

In Fermi coordinates the second fundamental form of $r$-level sets may be written as $(h_r)_{\beta \gamma} = \frac{1}{2} \partial_r g_{\beta \gamma}$.  We raise an index to obtain a family of shape operators $S(r)$, where $S(r)^{\beta}_{\gamma} = g^{\beta \nu} (h_r)_{\nu \gamma}$.  A computation shows that $S$ satisfies a Riccati equation involving curvature, namely
\begin{eqnarray}
\label{eqn:ricatti}
( \partial_r S(r) + S(r)^2 )_{\beta}^{\nu} = -R_{0 \beta \; 0}^{\; \; \; \; \nu}.
\end{eqnarray}

We will make use of Jacobi fields suitable to our coordinates.  Fix $p \in Y$ and consider the outward normal geodesic $\gamma_p$.  Choose any curve $\sigma$ in $Y$ such that $\sigma(0) = p$ and define a variation through geodesics by $\Gamma(s,t) = \E(\sigma(s), t )$.  This gives rise to a Jacobi field $J(t) = \partial_s \Gamma(s,t) |_{s=0}$ along $\gamma_p$.  Explicitly,
\[ J(t) = \dot\sigma^{\beta}(0) \partial_{\beta} |_{(\sigma(0),t)}. \]
So these special Jacobi fields have constant components in Fermi coordinates.  Convexity of $Y$ easily implies the following estimates:

\begin{lemma}  Let $J(t) = \dot\sigma^{\beta}(0) \partial_{\beta} |_{(\sigma(0),t)}$ be the Jacobi field along $\gamma$ described above.  Then\begin{enumerate}
\item $\inn{J(0)}{D_t J(0)} \geq 0$,
\item $|J(t)| \geq |J(0)| \; \cosh(at)$.
\end{enumerate} 
\label{lemma:jacobigrowth}
\end{lemma}

\begin{comment}
\begin{proof}
\begin{enumerate}
\item Observe that the symmetry lemma (see \cite[Lemma 6.3]{Lee}) and convexity imply
\begin{eqnarray*}
\inn{J(0)}{ D_t J(0)} & = & \inn{\dot\sigma(0)}{D_t \partial_s \Gamma(s,t)|_{(0,0)}} \\
                 & = & \inn{\dot\sigma(0)}{D_s \partial_t \Gamma(s,t)|_{(0,0)}} \\
                 & = & \dot\sigma^{\beta}\dot\sigma^{\nu} h_{\beta \nu} \\
                 & \geq & 0.
\end{eqnarray*}
\item The proof of the Jacobi field comparison theorem of \cite[Theorem 11.2]{Lee} is easily modified in the case when $J(0) \neq 0$ to give:
\[ |J(t)| \geq |J(0)| \; \cosh(at) + \frac{1}{a} \left(\frac{d}{dt} |J(t)|\right)_{t=0} \sinh(at). \]
Observe that part (1) now implies that $\left(\frac{d}{dt} |J(t)|\right)_{t=0} \geq 0$.
\end{enumerate}

Set $f(t) = |J(t)|^2$.  Then $f'(t) = 2 \inn{ D_t J(t)}{ J(t) }$ and 
\begin{eqnarray*}
f''(t) & = & 2 |D_t J(t)|^2 + 2 \inn{D^2_t J(t)}{ J(t) } \\ 
       & = & 2 |D_t J(t)|^2 + 2 \inn{ Rm( \sigma'(t), J(t)) \sigma'(t)}{ J(t) } \\
       & = & -2 sec( J(t), \sigma'(t) ) + 2 |D_t J(t)|^2.
\end{eqnarray*}
Since $f''(t) \geq 0$ by the curvature hypothesis, $f'$ is an increasing function.  As $f'(0) \geq 0$, $f'(t) \geq 0$ and so $f$ is increasing.  Therefore $f(t) \geq f(0)$.
\item 
\end{enumerate}
\end{proof}
\end{comment}

The comparison theorems we use are based on the treatment given in \cite{Petersen}.  These are obtained by analysis of the Riccati differential equation \eqref{eqn:ricatti}.  In what follows an inequality involving the shape operator of the form $S \geq c$  means that every eigenvalue of $S$ is greater than or equal to $c$.  Inequalities involving a metric are to be interpreted as inequalities between quadratic forms.

For the metric comparisons that follow we require a covering of the compact hypersurface $Y$.  Fix $\epsilon = \frac{1}{2} \min\{ \inj(Y), \pi \}$, where $\inj(Y)$ is the injectivity radius of $g_Y(y,0)$.  For any $y \in Y$, the ball $B^Y_{\epsilon}(y)$ is the domain of a convex normal coordinate chart with image $B_{\epsilon}(0) \subset \bR^n$.  On the ball $B_{\epsilon}(0)$, we will need to consider three metrics, the original $g_Y$ (transfered to $B_{\epsilon}(0)$ by means of normal coordinates), the round metric on the unit sphere $\bS^{n}$ in normal coordinates, and the flat metric in coordinates.  We will denote the round metric hereafter by $\grd$ and the flat metric by $\gbar$.  On compact subsets of $B_{\epsilon}(0)$ all three of these metrics are comparable.  Since $g_y(0,0)_{\beta \nu} = \grd(0)_{\beta \nu} = \gbar(0)_{\beta \nu} = \delta_{\beta \nu}$, continuity of the metrics implies we may find an $r = r(y)$ with $0 < r < \epsilon/2$ so that
\begin{itemize}\setlength{\itemsep}{10pt}
\item $\frac{1}{4} \grd_{\beta \nu} \leq (g_Y)_{\beta \nu} \leq 4 \grd_{\beta \nu}$  on $B_r(0)$,
\item $\frac{1}{4} \gbar_{\beta \nu} \leq (g_Y)_{\beta \nu} \leq 4 \gbar_{\beta \nu}$ on $B_r(0)$, 
\item $B_r(0) \subset B_{2r}(0) \subset B_{\epsilon}(0)$.
\end{itemize}

Compactness of $Y$ yields a finite subcover of the balls $B_{r(y)}(y)$ that cover $Y$.  Label these finitely many balls $W_i$.  Label the balls with the same centres and radius $2r(y)$ as $V_i$ and observe $W_i \subset \overline{W_i} \subset V_i$.  We refer to the covering of $Y$ by $\{W_i\}$ as the \textit{reference covering for Y}.

The choice of this covering ensures that distances between points in $W_i \subset Y$ with respect to the metrics $(g_Y)_i$, $\grd_i$, and $\gbar_i$ are all comparable.  We refer to this property again as the \textit{distance comparison principle}. \label{dcp}

In the theorem that follows and throughout the note we take eigenvalues of the metric $g_Y$ with respect to the euclidean metric $\gbar_i$ in normal coordinates for the $V_i$.  We let $\Omega_i$ denote the maximum eigenvalue of $g_Y(y,0)$ in each $\overline{W_i}$, and then set $\Omega = \max_i \; \Omega_i$.  Similarly, let $\omega$ be the minimum eigenvalue of $g_Y(y,0)$ over the cover $\overline{W_i}$.  As this covering and constants will be used throughout the paper, we always use normal coordinates along $Y$ in any choice of Fermi coordinates that follows.

An adaptation of the comparison theorems in \cite{Petersen} yields the following theorem.

\begin{theorem}[Comparison theorems]  Let $(y^{\beta}, r)$ be Fermi coordinates for $Y$ on $W_i \times [0,\infty)$ for an open set $W_i$ in the reference covering described above.  Let $\Lambda, \lambda$ denote the maximal and minimal eigenvalues of the shape operator over $Y$.  We have:
\label{thm:comparison}

\flushleft\textbf{ Shape operator estimate: } 
\begin{eqnarray}
\label{eqn:shapeinequality}
 a \tanh( a( r + L_1) ) \; \delta^{\beta}_{\gamma} \; \leq \; S^{\beta}_{\gamma}(y, r) \; \leq \; b \coth( b( r + L_2) ) \; \delta^{\beta}_{\gamma},
\end{eqnarray}
where $L_1 = \frac{1}{a} \tanh^{-1}(\frac{\lambda'}{a})$, $L_2 = \frac{1}{b} \coth^{-1}(\frac{\Lambda'}{b})$, and
\begin{eqnarray}
\label{eqn:constants}
\Lambda' = \begin{cases}  \Lambda  & \text{ if } \Lambda > b, \\
         2b & \text{ if } \Lambda \leq b . \end{cases} \nonumber \\
\lambda' = \begin{cases}  \lambda  & \text{ if } \lambda < a, \\
         \frac{a}{2} & \text{ if } \lambda \geq a . \end{cases}
\end{eqnarray}
\textbf{ Metric estimate: }
\begin{eqnarray}
\label{eqn:metric}
L_3 \; \cosh^2( a( r+ L_1 ) ) \; \delta_{\beta \nu} \leq g_{\beta \nu}(y, r) \leq L_4 \; \sinh^2( b( r+ L_2 ) ) \; \delta_{\beta \nu}, 
\end{eqnarray}
where $\displaystyle L_3 = \left( \frac{\omega}{\cosh^2{ a L_1 } } \right)$, $\displaystyle L_4 = \left( \frac{\Omega}{\sinh^2{ b L_2 } } \right)$, and $\Omega, \omega$ are described above.

\label{lemma:metriccomp}
\end{theorem}

\section{Essential Subsets}
\label{sec:tconvex}

In this section we provide a sufficient condition for the submanifold $K \subset M$ to be an essential subset.  We assume that $K$ is a compact $(n+1)$-dimensional Riemannian submanifold with boundary, such that $Y := \partial K$ is a smooth hypersurface that is convex with respect to the outward pointing unit normal.  We discuss a condition that ensures that $E\colon Y \times [0, \infty) \rarrow \overline{\MlK}$ is a diffeomorphism.  This property allows us to relax the hypothesis that $M$ be simply connected in the Anderson-Schoen result; essential subsets replace the requirement that the map $exp_p\colon T_p M \rarrow M$ be a diffeomorphism which is ensured by the Cartan-Hadamard theorem.

The basic idea of an essential subset $K$ is to capture the topology of the manifold $M$ inside $K$ in a such way that the normal geodesic flow off of the boundary of $K$ is a diffeomorphism of the outward normal bundle $Y \times [0, \infty)$ onto $\overline{\MlK}$.  Some sort of hypothesis on the topology and geometry of $M \backslash K$ is necessary in order for the exponential map to be injective, as can be seen by considering the orbit space obtained by taking the upper half space model of the two-dimensional hyperbolic plane under the action of the discrete group of dilations $G = \{ 2^n: n \in \bZ \}$.  Given a sufficiently small ball centred at any point along the `waist' of this quotient space, i.e. the image of $x=0$, one finds normal geodesics off the ball that intersect.

A subset $K \subset M$ is \textit{totally convex in M} if whenever $p, q \in K$ and $\sigma\colon[0,1] \longrightarrow M$ is a geodesic such that $\sigma(0)=p$, $\sigma(1) = q$ then $\sigma([0,1]) \subset K$.  The inclusion of $K$ into $M$ is a homotopy equivalence; see \cite{Klingenberg} for details.  It is interesting that totally convex sets play an important and somewhat analogous role in the of theory of souls of positively curved manifolds.  We again refer the interested reader to \cite{Klingenberg} and the references therein.

We have the following sufficient condition for an essential subset.

\begin{theorem} Let $K \subset M$ be a compact Riemannian submanifold with hypersurface boundary $Y$.  Suppose that
\label{thm:characterizationES}
\begin{itemize}
\item K is totally convex in M, and
\item $\sec(\MlK) \leq 0$.
\end{itemize}
Then $K$ is an essential subset for $M$.
\end{theorem}

\begin{proof}
As $K$ is totally convex, it is also geodesically convex ( i.e. $K$ contains a geodesic segment between any two of its points ).  It is well known that $K$ geodesically convex implies that $Y= \partial K$ is convex.  

We now check that the image $E(Y\times [0, \infty))$ is a subset of $\overline{\MlK}$.  To see this, notice that any normal geodesic $\gamma_p$ that re-enters $K$ must lie entirely inside $K$ since $K$ is totally convex, but this violates the fact $\gamma_p$ is a geodesic ray with an outward pointing tangent vector at $p$.  

Next, Jacobi field estimates and the nonpositive curvature assumption on $\MlK$ imply that $E$ is a local diffeomorphism on $Y \times [0, \infty)$.   We need only argue that $E$ is bijective.  Surjectivity of $E$ onto $\MlK$ is easy to see: for any point $q \in \MlK$ there is a closest point $p \in K$ to $q$, and it is straightforward to argue that $\gamma_p(t_0) = q$ for some $t_0$.  In order to argue injectivity of $E$, let 
\begin{equation*}
\begin{split}
\epsilon_0 = \sup \{ \epsilon > 0 \; | \; & E\colon Y \times [0, \epsilon) \rarrow \overline{\MlK} \; \text{is a diffeomorphism}  \\ 
 & \text{onto its image} \}. 
\end{split}
\end{equation*}
If $E$ is not injective then for $n > 0$ we have distinct points $p_n, q_n \in Y$ and $r$-values $t_n, s_n \in [0, \epsilon_0 + 1/n)$ such that $E(p_n, t_n) = E(q_n, s_n)$.  By the choice of $\epsilon_0$ it is impossible that both $t_n < \epsilon_0$ and $s_n < \epsilon_0$, so we may assume $ \epsilon_0 < s_n < \epsilon_0 + 1/n$.  It is straightforward to argue that a bound of the form $\epsilon_0 - 2/n < t_n < \epsilon_0 + 1/n$ holds as well.  By compactness of $Y$ we may pass to convergent subsequences and obtain points $p$ and $q$ such that $E(p, \epsilon_0) = E(q, \epsilon_0)$.  The points $p$ and $q$ are distinct as $p_n$ and $q_n$ are distinct and $E$ is a local diffeomorphism in a neighbourhood of $(p, \epsilon_0)$.  We set $m = E(p, \epsilon_0) = E(q, \epsilon_0)$.

We now show that $\gamma_p$ and $\gamma_q$ meet `head on', i.e. that $\gamma_p'(\epsilon_0) = -\gamma_q'(\epsilon_0)$.  If this is not the case, choose a vector $X \in T_{m} (\MlK)$ such that
\begin{equation}
\label{neginnerprod}
 g( X, \gamma_p'(\epsilon_0) ) < 0 \; \text{and} \; g( X, \gamma_q'(\epsilon_0) ) < 0.
\end{equation}

Consider any path $\sigma\colon (-\delta, \delta) \rightarrow \MlK$ for some $\delta > 0$ with $\sigma(0) = m$, $\sigma'(0) = X$.  Because $E$ is a local diffeomorphism, near $p$ we have a smooth curve $\sigma_Y^p\colon (-\delta, \delta) \rightarrow Y$ and a positive smooth function $r^p$ such that $\sigma(t) = E( \sigma_Y^p(t), r^p(t))$.  We may obtain a variation through $Y$-normal geodesics by considering the map $\Gamma(s,t) = E( \sigma_Y^p(s), t r^p(s) )$.  The first variation formula and equation \eqref{neginnerprod} above imply that $r^p(s)$ is a decreasing function of $s$; in particular for $s$ sufficiently small and positive, $r^p(s) < \epsilon_0$.  The same argument may be applied near $q$, and this implies for small enough $s$ that $E( \sigma_Y^p(s), r^p(s) ) = E( \sigma_Y^q(s), r^q(s) )$ which contradicts the choice of $\epsilon_0$.  

We have proved that $\gamma_p(\epsilon_0) = \gamma_q(\epsilon_0)$ and $\gamma_p'(\epsilon_0) = -\gamma_q'(\epsilon_0)$.  Therefore $\gamma_p$ is a geodesic such that $\gamma_p(0) = p \in K$ and $\gamma_p(2 \epsilon_0) = q \in K$, and so the image of $\gamma_p$ is contained in $K$ by total convexity, a contradiction.  Thus $E$ is a bijective local diffeomorphism, and consequently $E\colon [0,\infty) \rarrow \overline{\MlK}$ is a diffeomorphism.  This means $K$ is an essential subset of $M$.
\end{proof}

%%%%%%%%%%%%%%%%%%%%%%%%%%%%%%%%%%%%%%%%%%%%%%%%%%%%%%%%%%%%%%%%%%%%%%%%%%%%%%%%%%%%%%%%%%%%%%%%%

\section{The topological compactification}
\label{sec:tracey}

\label{sec:topcomp}

In this section we explain how to compactify $M$ given a choice of essential subset $K$ by extending the exponential map to take values in $\Mbar \backslash K = (\MlK) \cup \Minf$.  Fix an essential subset $K$ and define an extension $\Ebar\colon Y \times (0,1] \rarrow \Mbar \backslash K$ by

\begin{equation} \Ebar(p,s) = \begin{cases} \E(p, 2 \tanh^{-1}{s}) & \text{ if } s \in (0,1), \\
                                            [ \; \E( p, t): t \geq 0\; ] & \text{ if } s = 1. \end{cases}
\end{equation}\label{ebarref}

Recall that the notation $[\gamma]$ above represents the equivalence class of the geodesic ray $\gamma$ under the asymptotic equivalence relation, see page \pageref{page:asymp}.  We have also collapsed the normal component using a diffeomorphism.  We now verify that $\Ebar$ is a bijection.  Relative to the diffeomorphism $E\colon Y \times [0, \infty) \rarrow \MlK$ of the previous section, we write a generic curve $\sigma$ in $\MlK$ as $\sigma = (\sigma_Y, \sigma_r)$.

\begin{prop} $\Ebar$ is injective.
\label{prop:Ebaroneone}
\end{prop}
\begin{proof}
We must show that given distinct points $p, q \in Y$ that the normal geodesics $\gamma_p, \gamma_q$ have unbounded distance as a function of time.  It suffices to show given any $t > 0$ that the length of any curve $\sigma$ from $\gamma_p(t)$ to $\gamma_q(t)$ is bounded below by an unbounded function of time.

Suppose that $\sigma$ leaves the collar $Y \times [t/2, \infty)$.  Then the normal contribution of the length integral and the decomposition of the metric imply $len(\sigma) \geq t$.  In case that $\sigma$ remains in the collar, $len(\sigma) \geq len(\sigma^P)$, where $\sigma^P$ is the projection of $\sigma$ onto $Y \times \{t/2\}$, i.e. $\sigma^P(s) = (\sigma_Y(s), t/2)$.  Then Jacobi field estimates imply that 
$$len(\sigma) \geq len( \sigma^P ) \geq \cosh(a t/2) d_Y(p,q).$$ 
Therefore the length of any curve from $\gamma_p(t)$ to $\gamma_q(t)$ is bounded below by an unbounded function of time.\end{proof}

In order to show that $\Ebar$ is surjective we consider an untrapped geodesic ray $\sigma$ parametrized by arc length.  Recall that untrapped means that $\sigma$ eventually escapes every compact set.  Eventually $\sigma$ remains inside $\MlK$ and we take $\sigma(0)$ to be any point inside $\MlK$.  In this parametrization, the growth of $\sigma_r$ is bounded below by a linear function.

\begin{lemma} Let $\sigma = (\sigma_Y, \sigma_r)$ be an untrapped geodesic ray in $\MlK$ parametrized by arc length.  Then there exist constants $C, B, t_0 > 0$ such that
$$ \sigma_r(t)\geq C t + B, \mbox{ for all } t > t_0. $$
\label{lemma:geodesicrgrowth}
\end{lemma}
\begin{proof}
As a geodesic, the normal component of $\sigma$ satisfies
\[ \ddot\sigma_r + \dot\sigma^{\alpha} \dot\sigma^{\beta} \Gamma^0_{\alpha \beta} = 0, \]
where we have used the index conventions from Section \ref{sec:notation} and the decomposition of the metric $g = dr^2 + g_Y(r)$ to obtain the vanishing of the $\Gamma^0_{00}$ and $\Gamma^{0}_{\alpha 0}$-Christoffel symbols.  The $\Gamma^0_{\alpha \beta}$-Christoffel symbol is the scalar second fundamental form and so by our comparison results of Theorem \ref{lemma:metriccomp} we have:
\begin{equation}
\begin{split}
\label{eqn:sigmar}
\ddot\sigma_r & = - \dot\sigma^{\alpha} \dot\sigma^{\beta} \Gamma^0_{\alpha \beta}  \\
 & \geq a \tanh(a(r + L_1) ) |\dot\sigma_Y|^2_{g} \\
 & \geq 0. 
\end{split}
\end{equation}

Thus $\dot\sigma_r$ is nondecreasing.  Since $\sigma_r$ is eventually unbounded, there is a $t_0$ where $\dot\sigma_r(t_0) > 0$ and so $\dot\sigma_r(t) \geq \dot\sigma_r(t_0) > 0$ for $t > t_0$.  Upon integration we find that $\sigma_r(t) \geq \dot\sigma_r(t_0) (t - t_0) + \sigma_r(t_0)$.  
\end{proof}

We now find a candidate base point for a normal geodesic asymptotic to $\sigma$.
\begin{lemma} 
Let $\sigma_Y \colon [0, \infty) \rarrow Y$ be the projection of $\sigma$ onto $Y$.  Then $len( \sigma_Y ) < \infty$.
\end{lemma}
\begin{proof}
Since $\sigma$ is parametrized by arc length, $ \dot\sigma^{\alpha} \dot\sigma^{\beta} g_{\alpha \beta}(\sigma(t)) = |\dot\sigma_Y|_g \leq 1$.  The metric estimate of Theorem \ref{thm:comparison} implies
\begin{equation*}
L_3 \; \cosh^2( a( \sigma_r(t)+ L_1 ) ) \sum_{\alpha} (\dot\sigma_Y^{\alpha})^2 \leq \dot\sigma^{\alpha} \dot\sigma^{\beta} g_{\alpha \beta}( \sigma(t) ) \leq 1.
\end{equation*}
We may also consider the projected curve $\sigma_Y$.  Here the metric estimates imply
\begin{eqnarray*}
\dot\sigma^{\alpha} \dot\sigma^{\beta} g_{\alpha \beta}( \sigma_Y(t),0 ) \leq 4 L_4 \; \sinh^2( b L_2 )  \sum_{\alpha} (\dot\sigma_Y^{\alpha})^2.
\end{eqnarray*}
Combining these estimates and the result of Lemma \ref{lemma:geodesicrgrowth} we obtain:
\begin{eqnarray*}
\dot\sigma^{\alpha} \dot\sigma^{\beta} g_{\alpha \beta}( \sigma_Y(t) ) \leq  \frac{C(L_2, L_3, L_4)}{ \cosh^2( a(Ct + B) + L_1)}
\end{eqnarray*}
Integrating the square root of both sides we find that $len( \sigma_Y ) < \infty$.
\end{proof}

Since $len( \sigma_Y ) < \infty$, the completeness of $M$ implies that $\sigma_Y$ has a limit, and since $Y$ is closed this limit is a point $p \in Y$.  Let $\gamma_p(t)$ denote the outward normal geodesic emanating from $Y$ at $p$.  We now show that $\Ebar$ is surjective by showing that $\gamma_p$ is asymptotic to $\sigma$.

\begin{prop} $\Ebar$ is surjective.
\label{prop:Ebaronto}
\end{prop}
\begin{proof}
Given the untrapped geodesic ray $\sigma$ above, the previous lemma establishes the existence of a candidate normal geodesic $\gamma_p$ to represent the equivalence class $[\sigma]$ in $\Minf$.  We prove that $d( \gamma_p(t), \sigma(t) )$ remains bounded as $t \rarrow \infty$.  The triangle inequality implies that it is sufficient to show separately that 
\[d\big( (p, \sigma_r(t) ), (\sigma_Y(t), \sigma_r(t) ) \big) \; \text{and} \; d \big( (p, t), (p, \sigma_r(t)) \big)\] 
remain bounded as functions of time.

\textit{Step 1.} \textit{ $d( (p, t), (p, \sigma_r(t)) )$ is bounded as $t \rarrow +\infty$:}   In this situation, $d( (p, t), (p, \sigma_r(t)) ) = |t - \sigma_r(t)|$, so we must show the quantity $|t - \sigma_r(t)|$ is bounded as $t \rarrow \infty$.  By the fundamental theorem of calculus this is equivalent to the statement $1 - \dot\sigma_r(t) \in L^1(t_0, \infty)$ for $t_0$ sufficiently large:
\[ |t - \sigma_r(t)| \leq \int^{t}_{t_0} 1 - \dot\sigma_r(s) ds + |t_0 - \sigma_r(t_0)|,\]
Since $\dot\sigma_r^2 + |\dot\sigma_Y|^2 = 1$, the integrability of $1 - \dot\sigma_r(t)$ is related to the integrability of  $|\dot\sigma_Y|^2$, for
\[ 1 - \dot\sigma_r(t) \leq (1 - \dot\sigma_r(t))(1 + \dot\sigma_r(t)) = |\dot\sigma_Y|^2,\]
for large $t$.  Just as in the proof of Lemma \ref{lemma:geodesicrgrowth}, the special form of the $r$-component of the geodesic equation in Fermi coordinates and the estimates of Theorem \ref{thm:comparison} imply that $\ddot\sigma_r \geq 0$ and 
\[|\dot\sigma_Y|^2 \leq (1/a) \coth(a(r + L_1) \ddot\sigma_r.\]
The fundamental theorem of calculus applied to $\ddot\sigma_r$ implies that  $\ddot\sigma_r \in L^1(t_0, \infty)$ and consequently since $\coth(a(r + L_1)$ is bounded, shows that $|\dot\sigma_Y|^2 \leq C \ddot\sigma_r$ for large enough $t$.  Thus $|\dot\sigma_Y|^2$ is integrable and the quantity $|t-\sigma_r(t)|$ remains bounded as $t \rarrow +\infty$.

\textit{Step 2.}  \textit{$d( (p, \sigma_r(t)), (\sigma_Y(t), \sigma_r(t) ))$ is bounded as $t \rarrow +\infty$:}   For each $t_0$ consider the curve \[\tau^{(t_0)}(s) = ( \sigma_Y(s), \sigma_r(t_0) ), \mbox{ on} \; [t_0, \infty). \]
Clearly $ d( ( \sigma_Y(t_0), \sigma_r(t_0) ), (p, \sigma_r(t_0)) ) \leq len( \tau^{(t_0)}) $, and so it suffices to show that $len( \tau^{(t_0)})$ is bounded above by a constant independent of $t_0$.  To this end we use Jacobi field estimates based at the $r$-level set of value $\sigma_r(t_0)$.

In particular consider the Jacobi field along $\gamma_{\sigma_Y(s)}(t)$ given in Fermi coordinates as the constant vector field $$J_s(t) = (\dot\sigma_Y(s), t) \in T_{(\sigma_Y(s), t)} \overline{\MlK}. $$ See Section \ref{sec:notation} for a description of these special Jacobi fields.  In order to make our application of Lemma \ref{lemma:jacobigrowth} transparent, rescale the time parameter by $\lambda = t - \sigma_r(t_0)$.  Then $J_s( 0 ) = \dot\tau^{(t_0)}(s)$, and for $s \in (t_0, \infty)$, $t \in (\sigma_r(t_0), \infty)$ and $\lambda \in (0, \infty)$, Lemma \ref{lemma:jacobigrowth} implies
\[ |J_s(\lambda)| \geq |J_s(0)| \cosh( a \lambda ). \]
Therefore we may write:
\[ 1 \geq |\dot\sigma_Y(s)|_{g_Y(\sigma_r(s))} = |J_s(\sigma_r(s))| \geq |\dot\tau^{(t_0)}(s)| \cosh( a ( \sigma_r(s) - \sigma_r(t_0))). \]
Consequently using this estimate and the estimate on $\sigma_r$ from Lemma \ref{lemma:geodesicrgrowth}, we find that
\[ |\dot\tau^{(t_0)}(s)| \leq \frac{1}{\cosh( a( C s + B - \sigma_r(t_0)) )}, \]
for constants, $C, B$.  Upon integration of this expression we find
\begin{eqnarray*}
len( \tau^{(t_0)} ) &\leq& \int^{\infty}_{t_0}  \frac{1}{\cosh( a( C s + B - \sigma_r(t_0)) )}\\
                    &\leq& C(a) \tanh^{-1} (e^{ a( Cs + B - \sigma_r(t_0) ) } ) \big|^{s=\infty}_{s=t_0}\\
                    &\leq& C(a) \frac{\pi}{2}.
\end{eqnarray*}

Thus $len( \tau^{(t_0)})$ is bounded independent of $t_0$.  This completes step 2.

\end{proof}

The proof of the above proposition can be extended to yield a stronger result that will be useful in proving that the topology on $\Mbar$ is well defined.  In the lemma below we consider two hypersurfaces $Y_1, Y_2$ where $Y_i$ is the boundary of an essential subset $K_i$.  We use the notation $\gamma_{p'}$ to denote the normal geodesic to $Y_1$ emanating from the point $p' \in Y_1$, and $\sigma_{q'}$ to denote the normal geodesic to $Y_2$ emanating from the point $q' \in Y_2$.

\begin{lemma} Let $K_1$ and $K_2$ be essential subsets of $M$ with $Y_i = \partial K_i$.  Given any point $q \in Y_2$ there exists an open neighbourhood $V_q \subset Y_2$ of $q$ and $B > 0$ such that for every $q' \in V_q$:
\[d(\sigma_{q'}(t), \gamma_{p'}(t)) \leq B, \mbox{for all } t \geq 0,\]
where $\gamma_{p'}$ is the unique normal geodesic ray emanating from $Y_1$ that is asymptotic to $\sigma_{q'}$.
\label{lemma:unifbdd}
\end{lemma}
\begin{proof}
Let the exponential map $\E\colon Y_1 \times [0, \infty) \rarrow \MlK_1$ be denoted $\E_1$.

Fix $q \in Y_2$ and $R>0$.  Now $\sigma_{q}$ is eventually outside every compact set, so there exists $T \geq 0$ such that $\sigma_{q}(t) \in \E_1(Y_1 \times [R, \infty))$ for $t\geq T$.  Further the $r$-component of this curve, $(\sigma_{q})_r$, is eventually strictly monotone increasing so we may increase $T$  if necessary to ensure that  $(\dot\sigma_{q})_r(t) > 0$ for $t \geq T$.  By continuity of the exponential map, there is a precompact open ball $V_q$ in $Y_2$ about $q$ such that for any $q' \in V_q$ we have both $\sigma_{q'}(t)  \subset  \E_1( Y_1 \times [R, \infty) )$ and $\dot\sigma_{q'}(t)  >  0$ for all $t \geq T$.  For each such $q'$, let $p'$ be the unique element of $Y_1$ such that $\gamma_{p'}$ is asymptotic to $\sigma_{q'}$.  By compactness, for any $q' \in \overline{V_q}$, $d(\sigma_{q'}(t), \gamma_{p'}(t))$ is uniformly bounded for $t \in [0, T]$.  We need only check that a uniform bound holds for the tails of the geodesics emanating from $V_q$.

In order to estimate $d(\sigma_{q'}(t), \gamma_{p'}(t))$ for $t\geq T$, proceed as in the proof of Proposition \ref{prop:Ebaronto}.  By continuity of the exponential map and by shrinking $V_q$ if necessary all constants may be chosen independently of $q'$.
\end{proof}

We may now prove the topological part of Theorem \ref{thm:main}:

\begin{theorem}\label{thmmarsh}
If $M$ is a complete, noncompact Riemannian manifold with an essential subset $K$ with pinched negative curvature as in (\ref{curvasmp}), then $M$ admits a geometric compactification as a topological manifold with boundary.
\end{theorem}
\begin{proof}
Throughout the proof we use the notation for $\gamma$ and $\sigma$ as described preceding Lemma \ref{lemma:unifbdd}.

Suppose $K_1, K_2$ are two essential subsets of $M$.  The propositions above imply that we get bijections $\Ebar_i\colon Y_i \times (0,1] \rarrow \Mbar \backslash K_i$.  Endow each subset $\Mbar \backslash K_i$ with the topology $\tau_i$ that makes $\Ebar_i$ a homeomorphism.  We now show that these topologies are equivalent.  Let $K$ be a compact set such that $K \supset K_1 \cup K_2$.  Consider the identity map from $(\Mbar \backslash K, \tau_1) \rarrow (\Mbar \backslash K, \tau_2)$, i.e. the composition $\psibar = \Ebar_2^{-1} \circ \Ebar_1$.  To show that the topology on $\Mbar \backslash K$ is independent of $K_i$, it suffices to show that $\psibar$ is a homeomorphism.  By the symmetric roles of the $K_i$, it suffices to prove that $\psibar$ is an open map. 

We already have that $\psi = \E_2^{-1} \circ \E_1$ is a diffeomorphism.  We need only check that open neighbourhoods in $\tau_1$ of points in $\Minf$ are taken to open neighbourhoods in $\tau_2$.  Choose a basis element of the form $\Ebar_1(U \times (c,1])$ where $U$ is open in $Y_1$.  For every $[\gamma_p] \in \Minf \cap \Ebar_1(U \times (c,1])$ where $\gamma_p$ is asymptotic to $\sigma_q$ we must find a neighbourhood $V_q \subset Y_2$ of $q$ and $d > 0$ so that $\Ebar_2( V_q \times (d, 1] ) \subset \Ebar_1(U \times (c,1])$.  It is sufficient to show that $\Ebar_2( V_q \times (d, 1) ) \subset \Ebar_1(U \times (c,1))$; equivalently we may show $E_2( V_q \times (d, \infty) ) \subset E_1(U \times (c,\infty))$ for some different constants $c, d$.

Set $W := \Ebar_1(U \times (c,1))$.  The tail of $\sigma_q$ is eventually in W; we may choose a $T>0$ so that $\sigma_q(t) \in W$ for $t\geq T$.  Since $\psi$ is a diffeomorphism, for each $t > T$ we may get a ball $V_q(t)$ about $q$ and $\epsilon(t) > 0$ such that $\Ebar_2(V_q(t) \times (t-\epsilon(t), t+\epsilon(t))) \subset W$.  We may assume that for $t_2 > t_1$ we have $V_q(t_2) \supset V_q(t_1)$, and that the radii of these balls are less than the injectivity radius of $Y_2$.  Now $\cap_{t > T} V_q(t)$ is either a ball or is the singleton set $\{q\}$.  In case the intersection is a ball, $V_q$, we have that $\Ebar_2( V_q \times (T, 1) ) \subset W$, which completes the argument.  Otherwise choose $q_n \rarrow q$ such that $q_n$ enters $W$ and eventually leaves it.  Let $p_n$ be the corresponding points on $Y_1$ so that $\sigma_{q_n}$ is asymptotic to $\gamma_{p_n}$.  By compactness of $Y_1$, we may pass to a convergent subsequence and assume that $p_n \rarrow p_0$.  But now the uniform bound of Lemma \ref{lemma:unifbdd} and continuity of the exponential map imply that
\[ d( \sigma_{q}(t), \sigma_{p_0}(t) ) = \lim_{n\rarrow\infty} d( \sigma_{q_n}(t), \gamma_{p_n}(t) ) \leq B. \]
This means that $\sigma_q$ is asymptotic to $\gamma_{p_0}$, and so by injectivity of $\Ebar_1$, $p = p_0$.  This implies $p_n \rarrow p$, and so $p_n$ is eventually inside $U$, a contradiction.  Thus  $\cap_{t > T} V_q(t)$ contains a ball.  Therefore the topology on $\Mbar \backslash K$ is well defined.
\end{proof}

%%%%%%%%%%%%%%%%%%%%%%%%%%%%%%%%%%%%%%%%%%%%%%%%%%%%%%%%%%%%%%%%%%%%%%%%%%

\section{Regularity of the compactification}
\label{sec:compreg}

The results of the previous section described how to compactify $M$ as a topological manifold with boundary given a specific choice of essential subset $K$, and that the topology of the compactification is independent of $K$.  In this section we lay the groundwork and prove that the compactification $\Mbar$ has a $C^{a/b}$ structure.  In order to do this we will first have to describe our explicit comparison with hyperbolic space and how this relates to Fermi coordinates.  Next, since the manifold $\MlK$ is not complete and we estimate distances in $\MlK$ as compared to hyperbolic space we explain how to refine the reference covering for $Y$.  Just as in hyperbolic space with a compact set $K$ removed, we have to check that for points $p$ and $q$ far enough from $K$ but with closest points $p'$ to $p$ and  $q'$ to $q$ on $K$ sufficiently close, the geodesic segment from $p$ to $q$ remains in $\MlK$.  Such a refined covering will be called a \textit{special covering for Y}.

Given these geometric preliminaries we define a bounded metric $d_K$ on $\MlK$.  Given two essential subsets $K_1$, $K_2$, each endowed with a special covering for $Y_i$, we establish a $C^{a/b}$ comparability estimate for distances in a subset of $M\backslash (K_1 \cup K_2)$.  Then in the proof of the main theorem we explain why the distance estimate yields a $C^{a/b}$ structure for $\Mbar$.

We now describe our comparison geometry and modification of the metric comparison described at the end of Section \ref{sec:notation}.  In particular, consider the reference covering of $Y$ by small normal coordinate balls $W_i \subset \overline{W_i} \subset V_i$ as described preceding Theorem \ref{thm:comparison}.  In each $W_i$ we may use the metric $\grd_i$ to obtain a hyperbolic metric of constant curvature $-\lambda^2$ given by
\[ (h_{\lambda})_i = dr^2 + \frac{\sinh^2(\lambda r)}{\lambda^2} \grd_i. \] \label{hmetric}
We will call these metrics \textit{hyperbolic comparison metrics}.  A little algebra applied to the metric estimates of Theorem \ref{thm:comparison} implies:

\begin{theorem}[Hyperbolic Metric comparison]  Consider Fermi coordinates $(y^{\beta}, r)$ for $Y$ on $W_i \times [0, \infty)$.  There exists an $R = R( \Lambda, \lambda, \Omega, \omega, a, b)$ independent of $i$ such that for every $r > R$:
\label{thm:hcomp}
\begin{eqnarray}
\frac{ \sinh^2( a( r - R ) )} {a^2} \; \grdi_{\beta \nu} \; \leq g_{\beta \nu}(y,r) \; \leq \;\frac{ \sinh^2( b( r+ R ) ) }{b^2}\; \grdi_{\beta \nu}.
\end{eqnarray}
In particular, for any points $p, q \in W_i \times [R, \infty)$ such that a geodesic segment\footnote{Recall from Section \ref{sec:notation} that we use \textit{geodesic segment} to mean a distance minimizing geodesic curve between two points.} from $p$ to $q$ lies entirely inside $W_i \times [R, \infty)$ then:
\begin{equation}
d_a(p, q) \leq d(p,q) \leq d_b(p,q),
\end{equation}
where $d_{\lambda}$ is the distance in the hyperbolic comparison metric on $W_i \times [R, \infty)$ of constant curvature $-\lambda^2$ described above.
\end{theorem}

We now provide an adaptation of the estimates used in \cite{Anderson-Schoen}.  We first begin with some estimates in the two-dimensional hyperbolic plane of curvature $-\lambda^2$, $H^2(-\lambda^2)$.  Let $p, q \in H^2(-\lambda^2)$, and take measurements from a third point $x \in H^2(-\lambda^2)$.  Suppose that $s=d_{\lambda}(p,x), t = d_{\lambda}(q, x)$, and let $\theta$ be the angle between the radial geodesic connecting $x$ to $p$ and the radial geodesic connecting $x$ to $q$.  The well known law of hyperbolic cosines \cite{Petersen} yields a formula involving the distance between $p$ and $q$ and these parameters:
\begin{equation} \label{eqn:hypdist}
\cosh(\lambda \; d_{\lambda}(p, q)) = \cosh(\lambda s) \cosh(\lambda t) - \sinh(\lambda s) \sinh(\lambda t) \cos(\theta).
\end{equation}
We use this formula throughout this section.  In the special case that $\theta = \pi/2$ we obtain the hyperbolic Pythagorean theorem,
\[ \cosh(\lambda \; d_{\lambda}(p, q)) = \cosh(\lambda s) \cosh(\lambda t). \]

Assume $t \geq s > 2R$.  We have:

\begin{lemma} In a two-dimensional hyperbolic plane, there exist positive constants $c_1, c_2, c_3 > 0$ depending on $\lambda$ so that the following estimates hold:
\label{lemma:hcomp1}
\begin{enumerate}
\item $d_{\lambda}(p,q) \leq \begin{cases} s + t + \frac{2}{\lambda} \log \theta + c_1 & \text{when } e^{\lambda s} \theta \geq 1, \\
                                           t - s + c_2                                 & \text{when } e^{\lambda s} \theta \leq 4. \end{cases}$
\item $d_{\lambda}(p,q) \geq s + t + \frac{2}{\lambda} \log \theta - c_3$.
\end{enumerate}
\end{lemma}

The above lemma is proved by straightforward estimation of \eqref{eqn:hypdist} and is essentially the form that Anderson-Schoen obtained in \cite{Anderson-Schoen}.

For the next estimate, let $p, q, x \in H^2(-\lambda^2)$ be as before except assume that $s=t$.  We need to estimate the distance from $x$ to the geodesic segment $\sigma$ from $p$ to $q$.  Again straightforward estimation yields:

\begin{lemma} In a two-dimensional hyperbolic plane, there exists positive constants $c_4, c_5$ such that 
\label{lemma:hcomp2}
\[d_{\lambda}(x, \sigma) \geq \begin{cases} - \frac{1}{\lambda} \log \theta - c_4 & \text{when }\; e^{\lambda s} \theta \geq 1, \\ 
                                            s - c_5 & \text{when }\; e^{\lambda s} \theta \leq 4. \end{cases} \]
\end{lemma}
\begin{comment}
\begin{proof}
Let $u = d_{\lambda}(x, \sigma)$.  Then the geodesic segment realizing this distance must intersect $\sigma$ orthogonally.  Let $d = d_{\lambda}(p,q)$.  The hyperbolic Pythagorean theorem then implies that
\[ \cosh(\lambda u) = \frac{\cosh(\lambda s)}{ \cosh(\lambda d/2) }. \]
This may be estimated by
\[ e^{\lambda u} \geq \cosh(\lambda u) \geq \frac{1}{2} \frac{e^{\lambda s}}{e^{\lambda d/2}}. \]
This easily implies
\begin{eqnarray*}
u &\geq s - \frac{d}{2} - \frac{1}{\lambda} \log(2).
\end{eqnarray*}
Applying Lemma \ref{lemma:hcomp1}, we find that if $e^{\lambda s} \theta \geq 1$, 
\[ u \geq  - \frac{1}{\lambda} \log \theta - \frac{c_1}{2} - \frac{1}{\lambda} \log(2). \]
If $e^{\lambda s} \theta \leq 4$ then we find
\[ u \geq  s - \frac{c_2}{2} - \frac{1}{\lambda} \log(2). \]
\end{proof}
\end{comment}

We now convert the above estimates in hyperbolic space into estimates suited to our Fermi coordinates.

\begin{lemma}[Extended Anderson-Schoen estimates]  Consider Fermi coordinates $(y^{\beta}, r)$ for $Y$ on $W_i \times [R, \infty)$\footnote{Recall the constant $R$ used here is the constant from Theorem \ref{thm:hcomp}, which is used throughout this section.}. Let $p,q \in W_i \times [R, \infty)$.  Then there exist positive constants $\{c_j\}_{j=1}^{8}$ depending only on $R$ and the reference covering such that:
\begin{enumerate}
\item $d_{b}(p,q) \leq \begin{cases} s_p + s_q + c_1 \log d_Y(p',q') + c_2 & \text{when } e^{b s_p} d_Y(p',q') \geq 2, \\
                                     s_q - s_p + c_3                       & \text{when } e^{b s_p} d_Y(p',q') \leq 2. \end{cases}$
\item $d_a(p, q) \geq s_p + s_q + c_4 \log d_Y(p',q') - c_5$,
\end{enumerate}
where $p = (p', s_p)$, $q=(q', s_q)$ in coordinates, and $s_q \geq s_p$.

Also if $s_p = s_q$ and $\sigma$ is a geodesic segment in the hyperbolic comparison metric (cf. page \pageref{hmetric}) from $p$ to $q$, then the minimum $r$-value of $\sigma$ satisfies:
\begin{enumerate}
\item[3.] $\sigma_{r_{min}} \geq \begin{cases} -c_6 \log d_Y(p',q') - c_7, &\text{when}\; e^{b s_p} d_Y(p',q') \geq 2, \\
                                              s - c_8, & \text{when}\;e^{b s_p} d_Y(p',q') \leq 2. \end{cases} $
\end{enumerate}
\label{lemma:as}
\end{lemma}
\textbf{Note:} In order to avoid a proliferation of constants we reuse the labels $c_1$ through $c_8$ above, hence these constants are not the same as the constants in Lemmas \ref{lemma:hcomp1} and \ref{lemma:hcomp2}.
\begin{proof}
The points $p$ and $q$ lie in exactly one coordinate 2-plane $\Pi$ perpendicular to $Y$.  Distances between $p$ and $q$ in the hyperbolic comparison metrics are realized by geodesics lying entirely in $\Pi$ and so we may use Lemmas \ref{lemma:hcomp1} and \ref{lemma:hcomp2} specific to two dimensions stated above in our metric comparisons.  Further by the choice of reference covering and the distance comparison principle (cf. page \pageref{dcp}) the distance $\theta$ along the unit sphere is comparable to distance along $Y$.  From the choice of reference covering it follows that if $e^{b s_p} d_Y(p',q') \geq 2$, then $e^{bs} \theta \geq 1$; similarly if $e^{b s_p} d_Y(p',q') \leq 2$, then $e^{bs} \theta \leq 4$.
\end{proof}

The manifold $\MlK$ is not complete.  Therefore we need to be careful when applying comparison geometry to estimate distances in $\MlK$, for geodesic segments could potentially leave the manifold $\MlK$ entirely.  Fortunately the curvature assumptions imply that at least for points far enough from $Y$ whose nearest points on $Y$ are close enough, geodesic segments remain in the domain of a Fermi chart.  We now explain how to obtain these special charts.  For $x \in Y, \mu > 0$ and $t_0 > 0$, we call
\begin{eqnarray*}
TC(x, \mu, t_0) = \{ (y^{\beta}, t) \in W_i \times [0, \infty): d_{Y}(x, y) \leq \mu \; \mbox{and } t \geq t_0 \},  
\end{eqnarray*}
a \textit{truncated cylinder} about $x$ in Fermi coordinates $(W_i\times [0,\infty), (y^{\beta},t))$.  

\begin{lemma}[Double Buffer]  Fix $x \in W_i \times \{0\}$ in the domain of a Fermi coordinate chart.  Then there exist positive constants $\epsilon$, $\delta$, $T_{OB}$, $T_{IB}$ depending on $x$ and the constant $R$ from Theorem \ref{thm:hcomp} such that if we define
\begin{eqnarray*} OB & = & TC( x, \epsilon+\delta, T_{OB}),\; \mbox{the ``outer'' buffer}, \\
IB & = & TC( x, \epsilon, T_{IB}),\; \mbox{the ``inner'' buffer}, \end{eqnarray*}
then if $p, q \in IB$, the $g$-geodesic segment from $p$ to $q$ remains entirely inside OB.
\label{lemma:db}
\end{lemma}
\begin{proof}
We will determine the above constants such that if $p, q \in IB$, then
\begin{enumerate}
\item There is a curve $\sigma$ from $p$ to $q$ with $\sigma \subset OB$ such that $len(\sigma) \leq d_b(p,q)$.
\item For any curve $\sigma'$ from $p$ to $q$ that escapes $OB$, $len(\sigma') > d_b(p,q)$.
\end{enumerate}
This implies that a geodesic segment between $p$ and $q$ lies in $OB$, and hence in the domain of a Fermi chart.  See Figure \ref{fig:double}.

\begin{figure} 
 \centering 
 \includegraphics[width=0.6\textwidth]{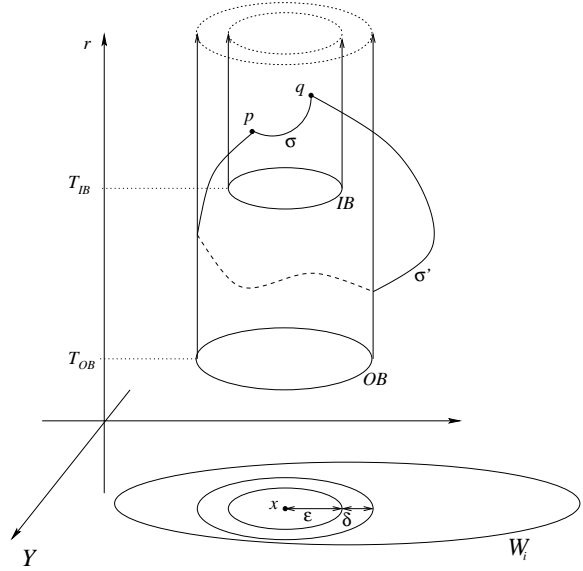} 
 \caption{Situation of Lemma \ref{lemma:db} }
 \label{fig:double} 
\end{figure}

\textit{Proof of Step 1:}  
Consider ``extremal'' choices of $p$ and $q$.  For $\epsilon$ and $T_{IB}$ to be determined, choose $p$ and $q$ to be any points on $B^Y_{\epsilon}(x) \times \{T_{IB}\}$ and write $p = (p', T_{IB})$, $q = (q', T_{IB})$.  Now set $\xi = d_Y(p',q')$.  Using Lemma \ref{lemma:as} we see that the $b$-hyperbolic geodesic segment $\sigma$ from $p$ to $q$ has minimal $r$-value $\sigma_{r_{min}}$ greater than or equal to
\[ \sigma_{r_{min}} \geq \begin{cases} -c_6 \log \xi - c_7, &\text{when}\; e^{b T_{IB}} \xi \geq 2 \\
                                       T_{IB} - c_8, & \text{when}\;e^{b T_{IB}} \xi \leq 2. \end{cases} \]
We impose the condition that $\sigma \subset OB$, i.e. that $\sigma_{r_{min}} \geq T_{OB}$.  Since $\xi \leq 2 \epsilon$,  this imposes the conditions:
\begin{equation}
\begin{split}
 \label{eqn:loccond1}
 T_{OB} & \leq  c_6 \log\left( \frac{1}{2\epsilon}\right) - c_7, \\
 T_{IB} & \geq  T_{OB} + c_8. 
\end{split}
\end{equation}

We now argue that there are no other conditions imposed on the constants after considering the extremal case.  Pick arbitrary points $u,v \in IB$ and let $\tau$ be the geodesic segment from $u$ to $v$ in the $b$-hyperbolic metric.  We must check that $\tau$ remains in $OB$.  Now $\tau$ must lie in a two-dimensional plane.  In this plane, let $p'$, $q'$ be the normal projections of $u, v$ on $ B^Y_{\epsilon}(x) \times \{0\}$, and consider $p=(p', T_{IB})$ and $q=(q',T_{IB})$ with $\sigma$ as before.  The two-dimensional region between the normal geodesics containing $u$ and $v$ and with $r$-values greater than those of $\sigma$ is convex in the hyperbolic metric.  This implies that $\tau$ remains in this region and so inside $OB$.

Thus given conditions \eqref{eqn:loccond1} above, the geodesic segment in the $b$-hyperbolic metric is a curve from $p$ to $q$ that remains in $OB$.

\textit{Proof of Step 2:}
Let $p, q \in IB.$  Suppose that $\sigma'$ is a curve between $p$ and $q$ that leaves $OB$.  The length of such a curve satisfies:
\[ len(\sigma')~\geq~d(p,\partial OB)~+~d(q,\partial OB).\]

The boundary of $OB$ is the union $\partial OB_1 \cup \partial OB_2$ where $\partial OB_1$ is the ``bottom'' disc $B^Y_{\epsilon+\delta}(x) \times \{T_{OB}\}$ and $\partial OB_2$ the ``vertical walls'', $\partial B^Y_{\epsilon+\delta}(x) \times [T_{OB}, \infty)$, for some $\delta$ to be determined.  Now $d(p, \partial OB_1) = r_p - T_{OB}$, as vertical geodesics are always minimizing.  

To estimate $d(p, \partial OB_2)$, suppose $\eta = (\eta', r_{\eta})$ is a point on $\partial OB_2$ closest to $p$.  Note that a $g$-geodesic segment from $p$ to $\eta$ must lie in $OB$, and so by Theorem~\ref{thm:hcomp} we may compare to the $a$-hyperbolic metric, i.e.
\[ d(p, \partial OB_2) \geq d_a(p, \partial OB_2). \]
Applying the comparison from Lemma \ref{lemma:as} we find
\begin{eqnarray*}
d_a(p, \partial OB_2) & \geq & r_p + r_{\eta} + c_4 \log{ d_Y(p', \eta') } - c_5 \\
                      & \geq & r_p + r_{\eta} + c_4 \log{\delta} - c_5 \\
                      & \geq & r_p + c_4 \log{\delta} - c_5,
\end{eqnarray*}
as $r_{\eta} \geq 0$.  Therefore:
\begin{eqnarray*}
d(p,\partial OB) \geq \min\left\{ r_p - T_{OB}, r_p + c_4 \log{ \delta } - c_5 \right\}.
\end{eqnarray*}
We impose the condition that $c_4 \log{1/\delta} + c_5 = -c_4 \log{\delta} + c_5 \leq T_{OB}$, and find that
\[ len(\sigma') \geq d(p,\partial OB) + d(q,\partial OB) \geq  r_p + r_q - 2 T_{OB}.\]
We need to apply Lemma \ref{lemma:as} once more to estimate $d_b(p,q)$.  Set $\xi = d_Y(p',q')$ and Lemma \ref{lemma:as} implies
\[ d_{b}(p,q) \leq \begin{cases} r_p + r_q + c_1 \log \xi + c_2 & \text{when } e^{b r_p} \xi \geq 2, \\
                                 r_q - r_p + c_3                & \text{when } e^{b r_p} \xi \leq 2. \end{cases}\]
Consider each case separately.  In order to guarantee that $len(\sigma') > d_b(p,q)$ we may therefore impose
\[ r_p + r_q - 2 T_{OB} > r_p + r_q + c_1 \log(\xi) + c_2,\]
in case $e^{b r_p} \xi \geq 2$, and 
\[ r_p + r_q - 2 T_{OB} > r_q - r_p + c_3,\]
when $e^{b r_p} \xi \leq 2$.  For the first case this is implied by the condition
\[ T_{OB} < \frac{c_1}{2} \log\left(\frac{1}{\xi}\right) - \frac{c_2}{2}, \]
and since this condition must hold for any $p, q \in IB$, it must hold when the $\log$-term is as small as possible, i.e. when $\xi = 2 \epsilon$.  So case 1 imposes 
\[ T_{OB} < \frac{c_1}{2} \log\left(\frac{1}{2 \epsilon}\right) - \frac{c_2}{2}. \]
The second case is implied by
\[ T_{OB} < r_p - \frac{1}{2} c_3 \leq \frac{1}{b}\log\left( \frac{2}{ \xi} \right) - \frac{1}{2} c_3,\]
and once again in order for this to hold for all $p, q$ we must impose:
\[ T_{OB} < r_p - \frac{1}{2} c_3 \leq \frac{1}{b}\log\left( \frac{2}{ 2 \epsilon} \right) - \frac{1}{2} c_3.\]

This completes the proof of Step 2.  To conclude the proof, we have seen that the conditions that must be satisfied in order to satisfy both items 1 and 2 above are:
\begin{enumerate}
\item $\delta + \epsilon$ is small enough so that $B^Y_{\delta+\epsilon}(x) \subset W_i$.
\item $T_{OB} \geq c_4 \log(\frac{1}{\delta}) + c_5$
\item $T_{OB} < \min \left\{ \frac{c_1}{2} \log\left(\frac{1}{2 \epsilon}\right) - \frac{c_2}{2},\frac{1}{b}\log\left( \frac{1}{ \epsilon} \right) - \frac{1}{2} c_3 , c_6 \log\left( \frac{1}{2\epsilon}\right) - c_7 \right\}.$
\item $T_{IB} \geq T_{OB} + c_8$. 
\end{enumerate}
Recall $\log(\frac{1}{z}) \rarrow +\infty$ as $z \rarrow 0^{+}$.  First choose $\epsilon$ and $\delta$ to meet condition 1; clearly any smaller $\epsilon$ will also work.  This fixes $\delta$, so now choose $T_{OB}$ subject to condition 2.  Shrinking $\epsilon$ if necessary, we may also satisfy condition 3. Finally choose $T_{IB}$ large enough to meet condition 4.
\end{proof}

Having finished the geometric preliminaries, we are now ready to describe the $C^{a/b}$ structure for $\Mbar$ that is independent of essential subset.  We begin by describing the basic philosophy of the proof.  In order to show that $\Mbar$ has a $C^{a/b}$ structure we must construct a $C^{a/b}$ atlas for $\Mbar$.  Given an essential subset $K_1 \subset M$, we use the double buffer lemma to obtain a collection of truncated cylinders that cover a neighbourhood of infinity in a sense that we make precise below.  We then obtain Fermi coordinates on these cylinders, and by collapsing the normal $r$-coordinate by a diffeomorphism, we obtain a coordinate cylinder that covers a deleted neighbourhood of the boundary $\Minf \subset \Mbar$.  We will show that transition functions from these cylinders to the collapsed truncated cylinders emanating from a second essential subset $K_2$ are $C^{a/b}$ functions.  As will be seen in the proof of Theorem \ref{thm:maincal} below, the transition functions will then extend by uniform continuity to $C^{a/b}$ functions on a coordinate cylinder including an open subset of $\Minf$.

Consider two essential subsets $K_1, K_2$ for $M$.  We begin with $K_1$.  By Theorem \ref{thmmarsh}, every point $p \in \Minf$ is the image under $\Ebar_1$ of exactly one point $p' \in Y_1$.  By Lemma \ref{lemma:db} we obtain parameters $\epsilon(p)$, $\delta(p)$, $T_{IB}(p)$, $T_{OB}(p)$.  Since the collection $\{B_{\epsilon(p)}(p')\}$ covers $Y_1$, we pass to a finite subcover 
\[ {\cB}_1 := \{B_{\epsilon(p_k)}(p_k')\}_{k=1}^{N_1}. \]
Set $T_1 = \max \{ T_{IB}(p_k): 0 \leq k \leq N_1 \}$.  Notice in Fermi coordinates relative to $Y_1$ if we write $p = E_1(p', r_p)$, $q= E_1(q', r_q)$ and assume $p',q' \in B_{\epsilon(p_k)}(p_k')$ for some $k$ and $\min\{r_p, r_q\} \geq T_1$ then a $g$-geodesic segment from $p$ to $q$ remains in some ``double buffer'' where we have comparison to hyperbolic metrics.  In what follows we only use this property and we will not mention the underlying double buffer structure explicitly again.

The same procedure may be repeated to obtain a collar neighbourhood of infinity relative to $Y_2$, and we let ${\cB}_2, N_2, T_2$ denote the corresponding data for $Y_2$ as described above for $Y_1$.  We set 
\begin{equation}
\label{Tdefn}
T = \max\left\{T_1 , T_2, \frac{1}{2a} ( 1 - \log(e-2) ) + \diam(K_1 \cup K_2)\right\}.
\end{equation}
The reason for the last term in the definition of $T$ will become apparent during the proof of Proposition \ref{prop:calphonm}.  We call $\cB_j$ the \textit{special coverings for $Y_j$}, $j=1, 2$, and the region $E_1( Y_1 \times (T, \infty) ) \cap E_2( Y_2 \times (T, \infty) )$ the \textit{special (deleted) neighbourhood of infinity}.  Observe also that every $p \in \Minf$ is in the intersection of the $\Mbar$-closure of two truncated cylinders, one emanating from each of $Y_1$ and $Y_2$.  As such, the truncated cylinders are deleted neighbourhoods of points in $\Minf$.  We introduce notation for these truncated cylinders and their images.  For $j = 1, 2$ let
\begin{align*}
TC^k_j & := B^k \times [T, \infty) \subset Y_j \times [0, \infty), \text{where } \cB_j = \{B^k\}_{k=1}^{N_j}, \\
\cTC_j & := \{ TC_j^k: 0 \leq k \leq N_j \}, \\
C_j^k & :=  E_j( TC_j^k ) \subset M, \\
\cC_j & := \{ C_j^k: 0 \leq k \leq N_j \}.
\end{align*}
Observe that in this notation, the lower index denotes the essential subset index and the upper index denotes an element of the special cover.

Since each $B \in \cB_j$ is contained in some $W_i$ from the reference covering for $Y_j$, $B$ is the image of a coordinate parametrization $\phi\colon \widetilde{B} \subset \bR^{n+1} \rarrow B \subset M$.  As $E\colon Y \times [0, \infty) \rarrow M$, we define $E_{coord} := E \circ (\phi \times Id)\colon \widetilde{B} \times [0, \infty) \subset \bR^{n+1} \rarrow M$.  

In what follows, we also consider the above constructions with $r$-coordinate collapsed by the diffeomorphism $\zeta\colon Y \times [0, \infty) \rarrow Y \times [0, 1)$ given by $\zeta(p, r) = (p, \tanh(r/2))$.  We use a circumflex to denote the collapsed version of subsets of $Y \times [0, \infty)$.  For example if $P\subset Y \times [0, \infty)$ then $\hat{P} = \zeta(P) \subset Y \times [0,1)$.  We also write the restriction of the map $\Ebar$ (cf. page \pageref{ebarref}) to $Y \times [0,1)$ as $\Ehat$.  Thus
\[ C_j^k = E_j( TC_j^k) = \Ehat_j( \widehat{TC}_j^k ). \]

To proceed we need to check that distances in the special neighbourhood of infinity measured relative to each essential subset are $C^{a/b}$ comparable.  This will be the key ingredient in showing that $\Mbar$ has a $C^{a/b}$ structure.  To facilitate this, given an essential subset $K_j$, for $p, q \in \MlK_j$ with $p = E_j(p', r_p)$, $q=E_j(q', r_q)$ in Fermi coordinates, define 
\begin{eqnarray}
d_{K_j} (p, q) = | e^{-r_p} - e^{-r_q} | + d_{Y_j} ( p', q' ).
\end{eqnarray}
This metric is defined on the entire set $\MlK_j$.  It is easy to verify that when restricted to a particular truncated cylinder $C^k_j$ with coordinate parametrization $(\Ehat_{coord})_j^k$, this metric is equivalent to the Euclidean metric in collapsed Fermi coordinates, i.e. that
\[ | e^{-r_p} - e^{-r_q} | + d_Y ( p', q' ), \mbox{and } \]
\begin{multline*}
| ((\Ehat_{coord})_j^k)^{-1}(p) - ((\Ehat_{coord})_j^k)^{-1}(q) | \\
= |( (p')^{\alpha}, \tanh(r_p/2)) - ((q')^{\alpha}, \tanh(r_q/2))|
\end{multline*}
are equivalent on $((\Ehat_{coord})_j^k)^{-1}(C^k_j)$.

We now prove the main $C^{a/b}$ distance estimate.
  
\begin{prop} There exists a positive constant $C$ depending on $a, b$ and $\diam(K_1 \cup K_2))$ such that whenever  $C_1^k \in \cC_1$ and $C_2^{k'} \in \cC_2$ and
\[ C_1^k \cap C_2^{k'} \neq \emptyset, \]
then
\begin{eqnarray}
\label{maincheck}
d_{K_2}( p, q ) \leq C \left(d_{K_1} (p, q) \right)^{a/b}, 
\end{eqnarray}
for all $p, q \in  C_1^k \cap C_2^{k'}$.
\label{prop:calphonm}
\end{prop}
\begin{proof}
Throughout the proof recall that we assume $a \leq 1 \leq b$.  We write $\alpha = a/b$.  

In Fermi coordinates relative to $K_1$ we write $p = E_1(p', r_p)$, $q=E_1(q',r_q)$, and with respect to $K_2$ we write  $p = E_2(\tilpp, \tilrp)$, $q=E_2(\tilqp,\tilrq)$.  By our assumption on $p, q$ and construction of the special covering and neighbourhood of infinity we have
\begin{equation}
\label{dcomp}
d_a(p, q) \leq d_M(p,q) \leq d_b(p, q),
\end{equation}
where $d_{\lambda}$ is the distance in the hyperbolic comparison metric (cf. page \pageref{hmetric}).    
By the distance comparison principle (cf. page \pageref{dcp}), $\theta$ is comparable to $d_{Y_1}( p',q')$, and $\tilth$ is comparable to $d_{Y_2}( \tilpp,\tilqp)$.  We are thus free to work with the angle $\theta$ and replace angles by a constant times distances along hypersurfaces upon obtaining the final estimate.

The inequality \eqref{dcomp} and the hyperbolic cosine law, equation \eqref{eqn:hypdist}, imply that
\begin{multline} 
\frac{1}{a} \cosh^{-1} \left(\cosh{(a\tilrp)\cosh(a\tilrq) } - \sinh{(a\tilrp)} \sinh{(a\tilrq)} \cos{\tilth} \right) \\
\leq \frac{1}{b} \cosh^{-1} \left(\cosh{ (b r_p) } \cosh{(b r_q)} - \sinh{(b r_p)} \sinh{(b r_q)} \cos{\theta} \right). 
\end{multline}

We use the estimate $1-\theta^2/2 \leq \cos(\theta) \leq 1-\theta^2/8$ for $0 \leq \theta \leq \pi$, and then the angle-sum formulas for hyperbolic cosine imply:
\begin{multline} 
\label{mainineq}
\frac{1}{a} \cosh^{-1} \left(\cosh{ (a(\tilrp - \tilrq)) } + \sinh{(a\tilrp)} \sinh{(a\tilrq)} \frac{\tilth^2}{8} \right) \\ 
\leq \frac{1}{b} \cosh^{-1} \left(\cosh{ (b(r_p-r_q)) } + \sinh{(b r_p)} \sinh{(b r_q)} \frac{\theta^2}{2} \right). 
\end{multline}
Set $D := \diam_M(K_1 \cup K_2)$. The triangle inequality implies that
\begin{equation} 
\begin{split}
\label{sstilcomp}
r_p - D & \leq \tilrp \leq r_p + D,\\
r_q - D & \leq \tilrq \leq r_q + D. 
\end{split}
\end{equation}

Assume that $r_q \geq r_p$ and the proof now breaks into two cases, first when $r_q-r_p \geq \log{(2)}$ and then when $0 \leq r_q - r_p \leq \log(2)$.

\subsubsection*{Case 1: $r_q-r_p \geq \log{(2)}$}

The main idea in this case is that we have $e^{-r_q} \leq \frac{1}{2} e^{-r_p}$, and therefore
\begin{eqnarray}
\label{case1main}
e^{-r_p} \leq 2 (e^{-r_p} - e^{-r_q}).
\end{eqnarray}
The main inequality \eqref{mainineq} above in conjunction with the estimate $$\alpha \cosh^{-1}(z) \leq \cosh^{-1}(z^{\alpha}),$$ valid for $z \geq 1$ and $0 \leq \alpha \leq 1$, imply
\begin{multline}
\label{case1red}
\cosh{ (a(\tilrp - \tilrq)) } + \sinh{(a\tilrp)} \sinh{(a\tilrq)} \frac{\tilth^2}{8} \\
 \leq  \left(\cosh{ (b(r_p-r_q)) } + \sinh{(b r_p)} \sinh{(b r_q)} \frac{\theta^2}{2} \right)^{\alpha}. 
\end{multline}
When $z \geq \frac{1}{2a} ( 1 - \log( e - 2 ) )$ we have 
\begin{equation} 
 e^{az - 1} \leq \sinh(az) \leq e^{a z}.
\label{tech:sinh}
\end{equation}
By our choice of $T$ in \eqref{Tdefn}, $r_q \geq r_p \geq \frac{1}{2a} ( 1 - \log( e - 2 ) )$.  Consequently:
\begin{align} 
\label{case1left}
\cosh{ (a(\tilrp - \tilrq)) } & + \sinh{(a\tilrp)} \sinh{(a\tilrq)} \frac{\tilth^2}{8}  \notag \\ 
&\geq \frac{e^{a(\tilrp-\tilrq)} + e^{a(\tilrq-\tilrp)}}{2} + e^{-2} e^{a(\tilrp+\tilrq)} \frac{\tilth^2}{8} \notag \\
&\geq \frac{e^{a(\tilrp+\tilrq)}}{8e^2} \left( e^{-2a\tilrq} + e^{-2a\tilrp} + \tilth^2 \right).
\end{align}
Similarly with the right hand side of inequality \eqref{case1red}, we may use the upper bound for hyperbolic sine provided by \eqref{tech:sinh} to obtain
\begin{align}
\label{case1right} 
\left(\vphantom{\frac{\theta^2}{2}}\cosh{ (b(r_p-r_q)) }\right. &\left. + \sinh{(b r_p)} \sinh{(b r_q)} \frac{\theta^2}{2}\; \right)^{\alpha} \notag \\ 
&\leq ( e^{b(r_p-r_q)} + e^{b(r_q-r_p)} + e^{b(r_p+r_q)}\theta^2)^{\alpha} \notag \\ 
&\leq   ( e^{b(r_p+r_q)} ( e^{-2br_q} + e^{-2br_p} + \theta^{2}) )^{\alpha}  \notag \\
&=  e^{a(r_p+r_q)} ( 2 e^{-2br_p} + \theta^{2} )^{\alpha}.
\end{align}
Combining inequalities \eqref{case1red}, \eqref{case1left}, \eqref{case1right}, dividing by $e^{a(\tilrp+\tilrq)}/(8e^2)$ and using \eqref{sstilcomp} to remove tildes gives:
\begin{eqnarray}
\label{ineq2}
e^{-2a\tilrq} + e^{-2a\tilrp} + \tilth^2 & \leq & 8e^{2aD+2} ( 2 e^{-2br_p} + \theta^{2} )^{\alpha}
\end{eqnarray}
An easy computation shows that we always have the estimate:
\begin{equation}
\label{tilest}
e^{-2 \tilrq} + e^{-2\tilrp} \geq ( e^{-\tilrq} - e^{-\tilrp} )^2.
\end{equation}
Recall that $a \leq 1$, and so $e^{-2az} \geq e^{-2z}$ for $z\geq 0$.  Apply this and inequality \eqref{tilest} to the left hand side of \eqref{ineq2} to obtain
\begin{eqnarray}
\label{ineq3}
e^{-2a\tilrq} + e^{-2a\tilrp} + \tilth^2 & \geq & e^{-2\tilrq} + e^{-2\tilrp} + \tilth^2 \notag \\
 & \geq & (e^{-\tilrp} - e^{-\tilrq})^2 + \tilth^2.
\end{eqnarray}
For the right hand side of \eqref{ineq2} we use $b\geq 1$ and the estimate \eqref{case1main} to see:
\begin{eqnarray}
\label{ineq4}
8 e^{2aD+2} ( 2 e^{-2br_p} + \theta^{2} )^{\alpha} & \leq & 8e^{2aD+2} ( 4(e^{-r_p} - e^{-r_q})^2 + \theta^2 )^{\alpha}.
\end{eqnarray}
Combining inequalities \eqref{ineq2}, \eqref{ineq3} and \eqref{ineq4} we have:
\begin{eqnarray}
 (e^{-\tilrp} - e^{-\tilrq})^2 + \tilth^2 & \leq & 8 \cdot 4^{\alpha} e^{2aD+2} ( (e^{-r_p} - e^{-r_q})^2 + \theta^2 )^{\alpha}.
\end{eqnarray}
This now implies \eqref{maincheck}, and completes the proof of Case 1.

\subsubsection*{Case 2: $0 \leq r_q-r_p \leq \log{(2)}$}

The main idea in this case is to use a power series expansion for hyperbolic cosine as $r_q-r_p$ is bounded.  Note that if $0 \leq r_q-r_p \leq \log{(2)}$ then 
\begin{equation} \label{eqn:rprq}
0 \leq |\tilrq - \tilrp| \leq \log(2) + 2D.
\end{equation}
Simple calculations imply that we may choose constants $k_1, \ldots, k_4$ so that for $0 \leq z \leq \log(2) + 2D$:
\begin{eqnarray}
1+k_1 z^2 \leq & \cosh{(z)} & \leq 1+ k_2 z^2, \label{locest1}\\
  k_3 z   \leq &  1-e^{-z}  & \leq k_4 z. \label{locest2}
\end{eqnarray}
So these estimates hold when $z = r_q - r_p$ or $z = |\tilrq-\tilrp|$.

We begin with inequality \eqref{case1red}.  We will first apply estimates \eqref{locest1} and the estimates for hyperbolic sine from \eqref{tech:sinh}.  We then apply the estimate $(1+x)^{\alpha} \leq 1 + x^{\alpha}$, valid for $x \geq 0$ and $0 \leq \alpha \leq 1$.  This yields
\begin{align}
1 + k_1 (a(\tilrp - \tilrq))^2 &+ \frac{1}{8e^{2}} e^{a(\tilrp+\tilrq)} \tilth^2 \nonumber \\
&\leq  1 + \left( k_2 (b(r_p-r_q))^2 + e^{b(r_p+r_q)}\theta^2 \right)^{\alpha}.
\end{align}
We cancel the ones and divide by $e^{a(\tilrp+\tilrq)}$, absorbing this factor into the right hand side of the inequality, obtaining
\begin{align}
\label{locest4}
k_1 e^{-a(\tilrp+\tilrq)}(a(\tilrp & - \tilrq))^2 + \frac{1}{8e^{2}} \tilth^2  \nonumber \\
& \leq  e^{2aD} \left( k_2 e^{-2br_p} (b(r_p-r_q))^2 + \theta^2 \right)^{\alpha}.
\end{align}
We consider the right hand side of this inequality.  A little algebraic manipulation, use of estimate \eqref{locest2}, and the fact that $b \geq 1$ give
\begin{align}
\label{locest5}
e^{2aD} ( k_2 e^{-2br_p} ( b (r_p & -r_q) )^2  + \theta^2 )^{\alpha} \notag \\
& \leq e^{2aD} \left( k_2 k_3^{-2} b^2 \left( e^{-r_p} ( 1 - e^{-(r_q - r_p)} ) \right)^2 + \theta^2 \right)^{\alpha} \notag \\
& \leq e^{2aD} \max\{(k_2 k_3^{-2} b^2)^{\alpha}, 1\} \left( (e^{-r_p}-e^{-r_q}) + \theta \right)^{2 \alpha}.
\end{align}
We now consider the left hand side of inequality \eqref{locest4}.  When $\tilrq > \tilrp$ we find:
\begin{align}
\label{locest6}
k_1 e^{-a(\tilrp  +\tilrq)} & (a(\tilrp  - \tilrq))^2 + \frac{1}{8e^{2}} \tilth^2  \nonumber \\
&\geq  a^2 k_1 k_4^{-2} e^{-2\tilrq} (1-e^{-(\tilrq-\tilrp)})^2 + \frac{1}{8e^{2}}\tilth^2 \nonumber \\
&\geq     a^2 k_1 k_4^{-2} (e^{\tilrp-\tilrq} (e^{-\tilrp}-e^{-\tilrq} ))^2 + \frac{1}{8e^{2}}\tilth^2 \nonumber \\
&\geq  \min\left\{a^2 k_1 k_4^{-2} e^{-2(\log(2) + 2D)}, \frac{1}{8e^{2}} \right\} \left( (e^{-\tilrp}-e^{-\tilrq})^2 + \tilth^2 \right),
\end{align}
where the last estimate uses inequality \eqref{eqn:rprq}.  We may now put estimates \eqref{locest5}, \eqref{locest6} together with \eqref{locest4} to get:
\begin{multline}
\label{locest7}
\min\left\{a^2 k_1 k_4^{-2} e^{-2(\log(2) + 2D)}, \frac{1}{8e^{2}} \right\} \left( (e^{-\tilrp}-e^{-\tilrq})^2 + \tilth^2 \right) \\ 
\leq e^{2aD} \max\{(k_2 k_3^{-2} b^2)^{\alpha}, 1\} \left( (e^{-r_p}-e^{-r_q}) + \theta \right)^{2 \alpha}.
\end{multline}
Observe that for positive $x, y$ and $z$, $x^2 + y^2 \leq z^2$ imply $x + y \leq 2z$.  Consequently this estimate and inequality \eqref{locest7} implies that for some constant $C$:
\begin{eqnarray}
e^{-\tilrp}-e^{-\tilrq} + \tilth & \leq & C \left((e^{-r_p}-e^{-r_q} ) + \theta \right)^{\alpha}.
\end{eqnarray}
An entirely similar computation holds for $\tilrp \geq \tilrq$.  This completes the proof of Case 2.

\end{proof}

We are now ready to finish the proof of Theorem \ref{thm:main}, which we obtain immediately from Theorems \ref{thm:characterizationES}, \ref{thmmarsh}, \ref{thm:hcomp} and the following:

\begin{theorem}\label{thm:maincal} Let $M$ be a complete Riemannian manifold containing an essential subset.  Suppose for every essential subset $K \subset M$ with reference covering $\{W_i\}$ for $Y=\partial K$ that there exists an $R>0$ such that for every $r>R$
\begin{eqnarray}
\frac{ \sinh^2( a( r - R ) )} {a^2} \; \grdi_{\beta \nu} \; \leq g_{\beta \nu}(y,r) \; \leq \;\frac{ \sinh^2( b( r+ R ) ) }{b^2}\; \grdi_{\beta \nu},
\end{eqnarray}
for all i, where $\grdi$ is the round metric in normal coordinates (cf. page \pageref{dcp}).  Then $\Mbar$ has a $C^{a/b}$ structure independent of $K$.
\end{theorem}
\begin{proof}
We must find an atlas of $C^{a/b}$ compatible charts.
  
Recall we have earlier defined the special neighbourhood of infinity $S = E_1( Y_1 \times (T, \infty) ) \cap E_2( Y_2 \times (T, \infty) )$, relative to essential subsets $K_1$ and $K_2$.  The complement $K_0 = M \backslash S$ is a compact set, and we choose an atlas ${\cC}_0$ of normal coordinate balls covering $K_0$ so that the collection of balls of half the radius still cover $K_0$.  Preceding Proposition \ref{prop:calphonm} we defined a covering $\cC_j, j= 1,2$.  Every truncated cylinder $C^k_j \in \cC_j$ is a deleted neighbourhood of points on the boundary of $\Minf$.  Let $\overline{C}^k_j$ be the union of $C^k_j$ and points of $\Minf$ in the $\Mbar$-closure of $C^k_j$; this is an open subset of $\Mbar$ containing an open subset of $\Minf$.  Set $\overline{\cC}_j = \{ \overline{C}^k_j: 1\leq k \leq N_j \}$.  We now show that ${\cC} = {\cC}_0 \cup \overline{\cC}_1 \cup \overline{\cC}_2$ is a $C^{a/b}$ compatible atlas for $\Mbar$.

Whenever a chart from $\cC_0$ overlaps with a chart from any $\cC_j, j = 0, 1, 2$ the transition function is smooth, and therefore $C^{a/b}$.  Similarly, transition functions from two charts in a single $\cC_j$, $j=1, 2$ are $C^{a/b}$ functions.

We now consider the case that a chart $\overline{C}^k_1 \in \overline{\cC}_1$ meets a chart $\overline{C}^{k'}_2 \in \overline{\cC}_2$.  But this is exactly the situation of Proposition \ref{prop:calphonm}.  We have a $C^{\alpha}$ estimate of the form:
\[ d_{K_2}( p, q ) \leq C \left(d_{K_1} (p, q) \right)^{\alpha}, \]
for points  $p, q \in C^k_1 \cap C^{k'}_2$.

Since $d_{K_j}(p, q) = | e^{-r_p} - e^{-r_q} | + d_{Y_j} ( p', q' )$ is equivalent to the Euclidean distance $| (\Ehat_{coord})_j^{-1}(p) - (\Ehat_{coord})_j^{-1}(q)|$ on $(\Ehat_{coord})_j^{-1} (C^k_1 \cap C^{k'}_2)$, we have that the transition function 
\[ \psi = ((\Ehat_{coord})^{k'}_2)^{-1} \circ (\Ehat_{coord})^k_1\]
is a $C^{\alpha}$ map on $((\Ehat_{coord})_1^k)^{-1} (C^k_1 \cap C^{k'}_2)$.  

As in the proof of Theorem \ref{thmmarsh}, $\psi$ extends to a continuous map $\psibar$ on the closure of $((\Ehat_{coord})_1^k)^{-1} (C^k_1 \cap C^{k'}_2)$, and by the result above extends to a $C^{a/b}$ map on the closure as well.  Thus
\[ \psibar = ((\Ebar_{coord})^{k'}_2)^{-1} \circ (\Ebar_{coord})^k_1\]
is a $C^{a/b}$ map on $((\Ebar_{coord})_1^k)^{-1} (\overline{C}^k_1 \cap \overline{C}^{k'}_2)$.

In summary, we have shown that given any essential subset $K$ we may construct a smooth atlas for $\Mbar$, and that any two such atlases are $C^{a/b}$ compatible.
These atlases are contained in a maximal atlas, which is a $C^{a/b}$ structure for $\Mbar$ independent of essential subset.
\end{proof}

\end{document}